\documentclass[11pt]{amsart}
\usepackage{amsfonts}
\usepackage{amsmath}
\usepackage[most]{tcolorbox}   % 加载 breakable 库
\usepackage{caption}           % 提供 \captionof

\usepackage{amssymb,amsthm,amscd}
\usepackage[mathscr]{eucal}
\usepackage{graphicx}
\usepackage{float}
\usepackage{color}
\usepackage[all]{xy}

\usepackage{graphicx}
\newtheorem{theorem}{Theorem}[section]
\newtheorem{lemma}{Lemma}[section]

\theoremstyle{remark}
\newtheorem{remark}{Remark}[section]
\theoremstyle{plain}

\title[THE $\lambda>0$ CASE FOR $\Delta u=\lambda-c e^u$ TYPE EQUATION ON HIGH DIMENSIONAL SPHERES]{THE $\lambda>0$ CASE FOR $\Delta u=\lambda-c e^u$ TYPE EQUATION ON HIGH DIMENSIONAL SPHERES}
\author{Yuxuan Li$^{1,2}$}
\address{$^1$School of Mathematics and Physics,
Xi'an Jiaotong-Liverpool University,
Suzhou, China}
\address{$^2$Department of Mathematical Sciences,
University of Liverpool,
Liverpool, UK}
\email{yuxuan.li2502@student.xjtlu.edu.cn}

\author{Qiang Niu$^{1,}$*}
\address{$^1$School of Mathematics and Physics,
Xi'an Jiaotong-Liverpool University,
Suzhou, China}
\email{Qiang.Niu@xjtlu.edu.cn}

\author{Wubin Zhou$^1$}
\address{$^1$School of Mathematics and Physics,
Xi'an Jiaotong-Liverpool University,
Suzhou, China}
\email{Wubin.Zhou@xjtlu.edu.cn}

\thanks{*Corresponding author.}

\begin{document}

\begin{abstract}
Motivated by prescribed Chern scalar curvature problems in complex geometry,
we study the equation
$\Delta_{\mathbb S^{N-1}}V+\lambda(e^V-1)=0$ with $\lambda>0$ on spheres.
For $N\ge4$ and $\lambda>2(N-2)$ with $\lambda\ne2N$,
we prove the existence of a nonconstant smooth solution
that is axially symmetric and symmetric about the equator. 

\end{abstract}


\maketitle

\setcounter{tocdepth}{1}
\tableofcontents
\bigskip

\section{Introduction}
The following equation on closed manifold has abundant applications in different areas of mathematics, such as conformal geometry \cite{KAZDAN,ANGELLA,FUSI,Win,yu,LIZHOU,KAZDAN1,DING3,ZOU,HO,LEJMI}, and abelian Chern-Simons-Higgs models \cite{TARAN1,TARAN2,TARAN3,NOL}.
\begin{equation}\label{main}
\Delta u=\lambda-c e^u . 
\end{equation}

On 2-dimensional closed Riemannian surface, the equation arises in the prescribed Gaussian curvature problem, where the nonlinear term $e^u$ determines
a conformal change of the background metric. The case when $\lambda<0$ could be solved by the standard maximum principle and priori estimate technique; the case when $\lambda>0$ is also called the mean-field equation, it is more subtle for the reason that we can not use maximum principle, we refer to \cite{KAZDAN,KAZDAN1,DING3,TARAN1,TARAN2,TARAN3,NOL,DING1,DJADLI} for related celebrated researches by Moser-Trudinger inequality and topological tools.

On higher dimensional closed manifold case, the equation arises in the prescribed scalar curvature problems in complex geometry \cite{ANGELLA,FUSI,Win,yu,LIZHOU,ZOU,HO,LEJMI}, here $e^u$ also determines
a conformal change of the background metric. The case when $\lambda<0$ could also be solved by the standard maximum principle and priori estimate technique; but when $\lambda>0$, the energy functional is unbounded from below \cite{ZOU}, which means that we can not use the traditional variational method to find the minimizer of the energy functional, there are also some studies using the flow approach \cite{ZOU,HO,LEJMI}, but the results obtained under
additional assumptions do not resolve the general
positive case. In prescribed curvature problems, $c$ and $\lambda$ are generally smooth functions, but in this paper we consider the positive constants case, aiming to obtain new
existence results and explore further phenomena in this positive setting.

We consider (\ref{main}) on high dimensional sphere, without loss of generality (by setting $V=u+\log (c / \lambda)$), we could assume $c=\lambda$, we consider the axially symmetric solutions to transfer the PDE to ODE. Besides, to better apply energy method, we consider the solution which is symmetric about the equator:
\begin{equation}
V^{\prime \prime}+(N-2) \cot \theta V^{\prime}+\lambda\left(e^V-1\right)=0,
\end{equation}
with boundary condition $V(\theta)=V(\pi-\theta), V^{\prime}(0)=V^{\prime}(L)=0$ and $L:=\frac{\pi}{2}$. We denote $V_a$ as the solution of the equation with initial condition
$$
V_a(0)=a, \quad V_a^{\prime}(0)=0.
$$

\section{Geometric significance and analytical difficulties}
\label{subsec:geometric-background}

Prescribed curvature problems ask whether a given curvature
function can be realized by a conformal change of a fixed
background metric.
In Hermitian geometry, the Chern--Yamabe problem seeks
a metric of constant Chern scalar curvature within a
prescribed conformal class \cite{ANGELLA}.
This is a natural approach to constructing distinguished
metrics on complex manifolds, particularly in the
non-K\"ahler setting.
Allowing the target curvature to vary leads to the
prescribed Chern scalar curvature problem \cite{FUSI},
while the corresponding problem for Gauduchon connections
extends this framework to almost Hermitian geometry
\cite{LIZHOU}.

The exponential nonlinearity follows directly from
the conformal transformation law.
Let $\omega$ be a Hermitian metric on a compact complex
manifold of complex dimension $m$, and let
$\widetilde\omega=e^u\omega$.
With the analysts' sign convention, the Chern Laplacian
on functions can be written as
\[
\mathcal D_\omega u
=\Delta_g u-\langle\vartheta,du\rangle_g,
\qquad
\Delta_g=\operatorname{div}_g\nabla,
\]
where $g$ is the associated Riemannian metric and
$\vartheta$ is the Lee form, defined by
\[
d\omega^{m-1}=\vartheta\wedge\omega^{m-1}.
\]
The Chern scalar curvature transforms according to
\[
S^{\mathrm{Ch}}(\widetilde\omega)
=e^{-u}\left(
S^{\mathrm{Ch}}(\omega)-\frac m2\mathcal D_\omega u
\right).
\]
Thus prescribing $S^{\mathrm{Ch}}(\widetilde\omega)=K$
is equivalent to
\[
\Delta_g u-\langle\vartheta,du\rangle_g
=\lambda-c e^u,
\qquad
\lambda=\frac{2}{m}S^{\mathrm{Ch}}(\omega),
\quad
c=\frac{2}{m}K.
\]
The coefficients therefore correspond, up to normalization,
to the curvatures before and after the conformal change. Restricting these coefficients to constants still leaves
a substantial nonlinear problem, since their positivity
retains the principal analytical difficulties discussed below.

The first-order term records part of the Hermitian
geometry, but the principal nonlinear difficulties
remain even when it vanishes.
Indeed, for a balanced metric,
$d\omega^{m-1}=0$, so $\vartheta=0$ and the equation
reduces exactly to
\[
\Delta_g u=\lambda-c e^u.
\]
Consequently, studying the equation without the drift
term retains a geometrically meaningful case.
It also retains the difficulties associated with positive
curvature, exponential growth, and concentration.
New existence or qualitative results for this equation
therefore have independent significance, even though
their extension to a general Hermitian background
requires additional analysis.

The sign of the prescribed curvature changes the
maximum-principle structure of the problem.
For two solutions $u_1,u_2$, their difference
$w=u_1-u_2$ satisfies
\[
\mathcal D_\omega w+c b(x)w=0,
\qquad
b(x)=\int_0^1e^{u_2+s(u_1-u_2)}\,ds>0.
\]
When $c<0$, the zeroth-order coefficient has
the sign needed for the standard uniqueness argument.
When $c>0$, this argument no longer applies.
The nonpositive Gauduchon-degree case of the
Chern--Yamabe problem was solved in \cite{ANGELLA},
and related prescribed curvature results were obtained
under suitable sign assumptions in \cite{FUSI,LIZHOU}.
The positive case is substantially more difficult,
and nonuniqueness can occur \cite{ANGELLA}.

A further difficulty is the loss of the exponential
integrability supplied by the natural energy in
dimension two.
On a closed Riemannian surface, the Moser--Trudinger
inequality gives an estimate of the form
\[
\log\left(
\frac{1}{\operatorname{Vol}_g(M)}
\int_M e^{u-\bar u}\,dV_g
\right)
\leq C_1\int_M|\nabla u|^2\,dV_g+C_2,
\,
\bar u=\frac{1}{\operatorname{Vol}_g(M)}
\int_Mu\,dV_g.
\]
This controls the exponential term by the Dirichlet
energy and is a fundamental tool in the analysis of
mean-field and prescribed Gaussian curvature equations.
For real dimension $d>2$, however, the embedding
$$
H^1(M)\hookrightarrow L^{2d/(d-2)}(M)
$$
provides only power integrability.
There is no analogous uniform bound of the displayed
form using only the quadratic Dirichlet energy.
Higher-dimensional Moser--Trudinger inequalities involve
$W^{1,d}$ control, which is not supplied by this energy.
Thus the two-dimensional variational approach does not
extend directly to higher dimensions.

This obstruction is reflected in the positive
Chern--Yamabe problem:
in the balanced setting with positive total Chern
scalar curvature, the normalized energy functional
is unbounded from below when the complex dimension
is at least two \cite{ZOU}.
Direct global minimization is therefore unavailable.
Geometric flows offer another approach
\cite{ZOU,HO,LEJMI}, but require a priori estimates
strong enough to establish long-time existence and
convergence; energy monotonicity alone is insufficient.
The cited works obtain results under additional
hypotheses without resolving the general positive case.
Understanding this regime is consequently important
both for constructing conformal metrics with prescribed
curvature and for describing concentration, nonuniqueness,
and compactness in exponential elliptic equations.

In this paper we study the equation on the unit round sphere \(S^{N-1}\) with \(N\ge4\). Within this range, \(S^6\) is the only sphere admitting an almost complex structure. The general Chern scalar curvature equation includes a first-order drift term determined by the Lee form. Our aim is to investigate a basic model that isolates the exponential nonlinearity and the positive-parameter regime. We therefore omit the first-order drift term and focus on the resulting equation on the round sphere, where symmetry allows a detailed analysis of the existence and profiles of solutions.

A drift term-free spherical equation can nevertheless arise by restricting the conformal factor on a different Hermitian manifold. Consider the standard Hopf surface from (\cite{ANGELLA}, Example 5.1),
\[
H=(\mathbb C^2\setminus\{0\})/\langle z\mapsto2z\rangle,
\qquad
g_H=|z|^{-2}g_{\mathrm{Eucl}}=dt^2+g_{S^3},
\]
where \(t=\log|z|\) modulo \(\log2\), and \(g_{S^3}\) is the unit round metric. In the conventions used here, direct calculation gives
\[
S^{Ch}(g_H)=4,\qquad \vartheta_H=-2dt.
\]
For a conformal factor \(V\) depending only on the \(S^3\) variable, the Lee term vanishes and
\[
D_{g_H}V=\Delta_{S^3}V,\qquad
S^{Ch}(e^Vg_H)=e^{-V}\bigl(4-\Delta_{S^3}V\bigr).
\]
Consequently, prescribing \(S^{Ch}(e^Vg_H)\equiv4\) is equivalent to
\[
\Delta_{S^3}V+4(e^V-1)=0.
\]
This example gives a geometric origin for a spherical equation without a drift term, even though the Hermitian background itself is not balanced. The conformally deformed metrics are defined on \(H\), with conformal factors depending only on its spherical variable. Motivated by such curvature equations, we study the positive-parameter problem on \(S^{N-1}\) for general \(N\ge4\), which is also of independent interest in nonlinear elliptic PDEs.

\section{Numerical implications and main Theorem}
For each fixed $\lambda>0$, we regard the north-pole value
$a=V(0)$ as the shooting parameter. More precisely, let $V_a$
denote the regular solution of
\[
V_a''+(N-2)\cot\theta\,V_a'
+\lambda(e^{V_a}-1)=0,
\qquad
V_a(0)=a,\quad V_a'(0)=0.
\]
The energy estimate will ensure that $V_a$ exists up to $L$.
Thus the problem is to find a nonzero initial value $a$
for which
\[
V_a'(L)=0.
\]
The choice $a=0$ gives only the trivial solution.

In the numerical experiments, we therefore vary the pole value $a$
and solve for $\lambda$ subject to $V_a'(L)=0$.
For bounded $\lambda$, the mass identity
$\int_{\mathbb S^{N-1}}e^V\,d\sigma=|\mathbb S^{N-1}|$,
together with the radial flux estimate, excludes positive blow-up
away from the poles, since such blow-up would force the exponential
mass on a fixed latitude band to diverge.

We use the damped Newton method to draw the graph, we pin $u(0)=a$ and
treat $\lambda$ as an additional unknown.

Inspired by the numerical experiments, we consider positive and negative shooting parameters separately.
For $2(N-2)<\lambda<2N$, we use $a>0$ to construct a solution
strictly decreasing from the north pole to the equator.
For $\lambda>2N$, we use $a<0$ to construct a solution
strictly increasing on the same interval.
\begin{figure}[htbp]
    \centering
    \includegraphics[width=0.7\textwidth]{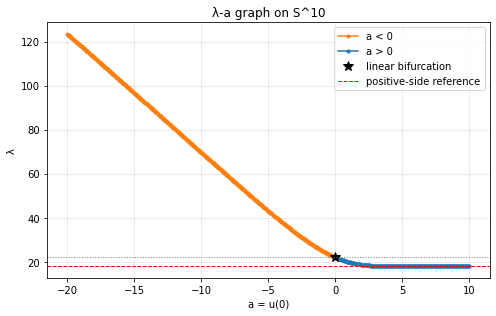}
    \caption{Bifurcation diagrams for $\Delta u=\lambda(1-e^u)$ on $S^{10}$.}
    \label{fig:gelfand}
\end{figure}

\begin{figure}[htbp]
    \centering
    \includegraphics[width=0.7\textwidth]{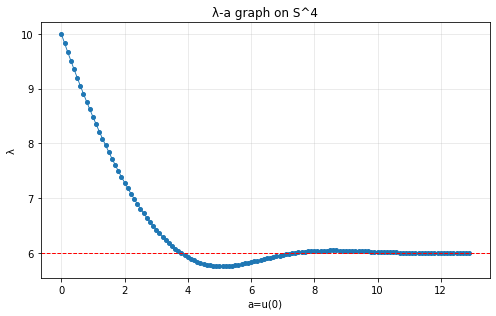}
    \caption{Bifurcation diagrams for $\Delta u=\lambda(1-e^u)$ on $S^{4}$.}
    \label{fig:gelfand}
\end{figure}

\begin{figure}[htbp]
    \centering
    \includegraphics[width=0.7\textwidth]{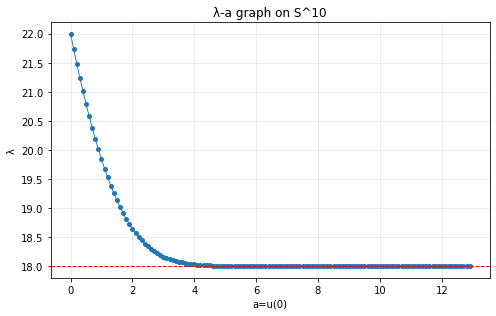}
    \caption{Bifurcation diagrams for $\Delta u=\lambda(1-e^u)$ on $S^{10}$.}
    \label{fig:gelfand}
\end{figure}

The numerical diagrams reveal different behaviors for positive
and negative pole values $a=V(0)$.
For $a>0$, the computed branch on $\mathbb{S}^{4}$ exhibits
damped oscillations about $\lambda=2(N-2)$, whereas no such
oscillations are observed on $\mathbb{S}^{10}$ ($N=11$)
over the computed range.
This contrast may be related to the dimensional threshold
in the classical Gelfand problem
$\Delta u+\lambda e^u=0$ with homogeneous Dirichlet boundary
conditions in the unit ball of $\mathbb{R}^{N}$ (bifurcation diagram is shown in Figure 4), whose radial
solution branch has infinitely many turning points for
$3\leq N\leq9$ but none for $N\geq10$
\cite{JOSEPH,DUPA}. The crucial difference between our problem and the classical Gelfand problem is that by maximum principal, we could only pin $u(0)=a>0$ in $\Delta u+\lambda e^u=0$ with homogeneous Dirichlet boundary
conditions, however, no sign restriction on $u(0)$ is imposed on the solutions in our problem, so we can study the case when $a<0.$

\begin{figure}[htbp]
    \centering
    \includegraphics[width=0.9\textwidth]{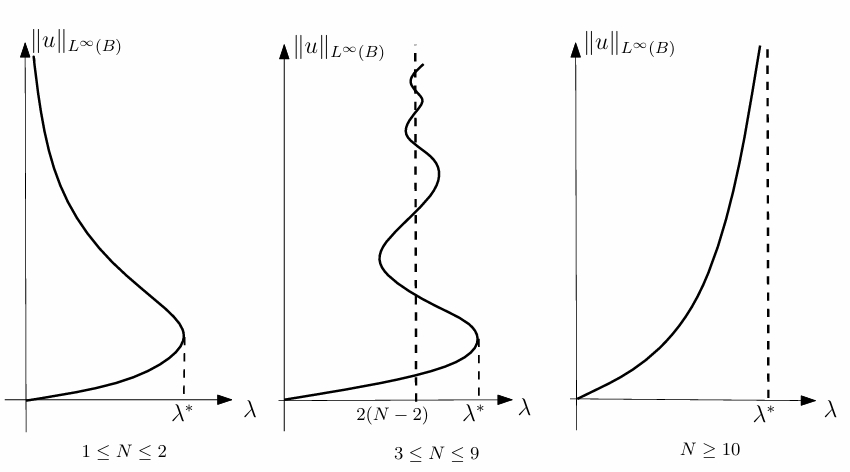}
    \caption{Bifurcation diagrams for $\Delta u+\lambda e^u=0$ with homogeneous Dirichlet boundary
conditions in the unit ball of $\mathbb{R}^{N}$.}
    \label{fig:gelfand}
\end{figure}

In Figure 1-3, for $a>0$, the computed positive-$a$ branch remains
below $\lambda=2N$, approaching this value as $a\downarrow0$.
For $a<0$, the numerical branch instead extends from
$\lambda=2N$ toward larger values of $\lambda$ as $a$ decreases.
Although the oscillations in lower dimensions allow the
positive-$a$ branch to enter the region $\lambda<2(N-2)$,
the value $2(N-2)$ remains theoretically distinguished,
as the parameter associated with the singular solution
$V(\theta)=-2\log\sin \theta+\log \frac{N-3}{N-2}$, besides, in \cite{miyamoto}, it is proved that, for \(4\le N\le10\), infinitely many nonconstant smooth axially symmetric solutions exist at \(\lambda=2(N-2)\), which can be explained by the infinite winding of the associated shooting orbit around the equilibrium value $2(N-2)$.
For positive shooting parameters, the present paper focuses
on proving existence in the interval $2(N-2)<\lambda<2N$.
\begin{theorem}\label{thm:nonconstant-symmetric-solutions}
Let $N\ge4$ and
\[
\lambda\in\bigl(2(N-2),\infty\bigr)\setminus\{2N\}.
\]
Then the equation
\[
\Delta_{\mathbb S^{N-1}}V+\lambda(e^V-1)=0
\quad\text{on }\mathbb S^{N-1}
\]
admits a nonconstant smooth solution which is axially symmetric
and symmetric about the equator.

Moreover, the solution can be chosen so that:
\begin{itemize}
    \item if $2(N-2)<\lambda<2N$, then
    \[
    V(0)>0,\qquad
    V'(\theta)<0
    \quad\text{for }0<\theta<\frac{\pi}{2};
    \]
    \item if $\lambda>2N$, then
    \[
    V(0)<0,\qquad
    V'(\theta)>0
    \quad\text{for }0<\theta<\frac{\pi}{2}.
    \]
\end{itemize}
\end{theorem}

\medskip
\noindent
\textit{\textbf{Outline of the proof of
Theorem~\ref{thm:nonconstant-symmetric-solutions}.}}
Set $L=\pi/2$, and let $V_a$ denote the regular solution
with $V_a(0)=a$ and $V_a'(0)=0$.
We seek a nonzero shooting parameter $a$ such that
$V_a'(L)=0$.

For $2(N-2)<\lambda<2N$, Lemma~4.1 shows that
$V_a'<0$ on $(0,L]$ for all sufficiently small positive $a$.
Lemma~4.2 provides a singular limit along a sequence
$a_j\to+\infty$, while Lemma~4.3 shows that this limit
turns upward before the equator.
Consequently, for sufficiently large $j$,
$V_{a_j}'$ is positive at an interior point.
Combining these facts with continuous dependence,
the shooting construction in Lemma~4.4 produces
$a_+>0$ such that
\[
V_{a_+}'(\theta)<0\quad(0<\theta<L),
\qquad V_{a_+}'(L)=0.
\]
Reflection across the equator, together with regularity
at the poles, gives the required smooth solution
on $\mathbb S^{N-1}$.

For $\lambda>2N$, the negative shooting construction
is established in
Theorem~\ref{thm:existence-above-2N}.
Sufficiently negative initial values yield
$V_a'>0$ on $(0,L]$, whereas linearization at the zero
solution shows that, for negative $a$ sufficiently close
to zero, $V_a'$ is negative at an interior point.
Taking the boundary of the interval of increasing
shooting solutions yields a solution with negative
pole value and strict increase on $(0,L)$.

Thus, Lemma~4.4, followed by reflection across the equator,
and Theorem~\ref{thm:existence-above-2N}
together prove
Theorem~\ref{thm:nonconstant-symmetric-solutions}.
\section{Pin initial value $V(0)=a>0$}

\medskip
\noindent\textit{Regularity at the poles.}
For any $a\in\mathbb R$, the integral formulation of the
shooting equation gives a unique local solution by a
contraction argument on a sufficiently small interval;
the Lipschitz constant of the integral operator is
$O(\delta^2)$ on a fixed neighborhood of the constant
function $a$.

Set
\[
b(\theta)=\left(\frac{\sin\theta}{\theta}\right)^{N-2},
\qquad b(0)=1,
\]
and let $W_a$ denote the even extension of $V_a$ near zero.
Changing variables in the integral formula gives, for
$\theta\ne0$,
\[
W_a'(\theta)
=-\lambda\theta\,b(\theta)^{-1}
\int_0^1 s^{N-2}b(s\theta)
\bigl(e^{W_a(s\theta)}-1\bigr)\,ds.
\]
The right-hand side extends continuously through zero
and is $O(\theta)$. Integrating this identity shows that
$W_a(\theta)-a=O(\theta^2)$, so $W_a$ is $C^1$ with
$W_a'(0)=0$.
If $W_a$ is $C^k$, differentiation under the integral sign
shows that the right-hand side is $C^k$, and hence $W_a$
is $C^{k+1}$. Induction therefore proves that $W_a$ is
smooth and even.

A smooth even function of geodesic distance defines a
smooth radial function in geodesic normal coordinates.
Thus the shooting solution is smooth at the north pole.
For the solution obtained by equatorial reflection,
smoothness at the south pole follows in the same way.
\medskip

\begin{lemma}
For sufficiently small positive $a$, we have $$
\lambda<2 N \quad \Longrightarrow \quad V_a^{\prime}(\theta)<0 \quad(0<\theta \leq L)
$$
\end{lemma}
\begin{proof}
Rewrite the original equation as
\begin{equation}\label{integral}
V_a^{\prime}(\theta)=-\lambda \sin ^{-(N-2)} \theta \int_0^\theta \sin ^{N-2} s\left(e^{V_a(s)}-1\right) \mathrm{d} s .
\end{equation}
Using the fact that $$
\cot \theta V_a^{\prime}(\theta) \rightarrow V_a^{\prime \prime}(0) ,
$$the expansion at north pole gives
$$
V_a^{\prime \prime}(0)=-\frac{\lambda}{N-1}\left(e^a-1\right)<0 .
$$
Define
$$
G(q):=e^q-q-1 \geq 0, \quad E_a(\theta):=\frac{1}{2}\left(V_a^{\prime}\right)^2+\lambda G\left(V_a\right) .
$$
Consider
$$
E_a^{\prime}(\theta)=-(N-2) \cot \theta\left(V_a^{\prime}\right)^2 \leq 0, \quad 0<\theta<L,
$$
we will have
$$
\frac{1}{2}\left(V_a^{\prime}\right)^2+\lambda G\left(V_a\right) \leq \lambda G(a) .
$$
Since $G(q) \rightarrow+\infty$ as $|q| \rightarrow \infty$, we will have the boundness of $V_a$ and $V_a^{\prime}$ on $[0, L]$.

Define
$$\Phi_\lambda(\theta)=\left.\frac{\partial V(\theta ; a)}{\partial a}\right|_{a=0},$$ then differentiate the original equation w.r.t $a$ gives
\begin{equation}\label{diffa}
\left(\partial_a V_a\right)^{\prime \prime}+(N-2) \cot \theta\left(\partial_a V_a\right)^{\prime}+\lambda e^{V_a} \partial_a V_a=0 ,
\end{equation}
if we denote $\Phi_\lambda=\left.\partial_a V_a\right|_{a=0}$ then (\ref{diffa}) becomes
\begin{equation}\label{equationa}
\left\{\begin{array}{l}
\Phi_\lambda^{\prime \prime}+(N-2) \cot \theta \Phi_\lambda^{\prime}+\lambda \Phi_\lambda=0, \\
\Phi_\lambda(0)=1, \\
\Phi_\lambda^{\prime}(0)=0 .
\end{array}\right.
\end{equation}
To justify this step, let \(\Phi_\lambda\) denote the regular solution of (\ref{equationa}). The energy estimate gives \(\|V_a\|_{C^1([0,L])}\le C|a|\) for sufficiently small \(|a|\), and hence \(e^{V_a}-1-V_a=O(a^2)\) uniformly on \([0,L]\). Subtracting the integral equations for \(V_a/a\) and \(\Phi_\lambda\) and applying Grönwall’s inequality yields \(\|V_a/a-\Phi_\lambda\|_{C^1([0,L])}\le C|a|\to0\). Since \(V_0\equiv0\), this establishes differentiability of the solution map at \(a=0\) and identifies \(\left.\partial_aV(\cdot;a)\right|_{a=0}=\Phi_\lambda\), justifying both the parameter derivative and the asserted convergence.

Differentiate (\ref{equationa}) and define $z_\lambda:=-\Phi_\lambda^{\prime}$, we get
\begin{equation}\label{equationz}
-\left(\sin ^{N-2} \theta z_\lambda^{\prime}\right)^{\prime}+(N-2) \sin ^{N-4} \theta z_\lambda=\lambda \sin ^{N-2} \theta z_\lambda .
\end{equation}
When $\lambda=2 N$, $z_*(\theta):=\sin \theta \cos \theta$ is a positive solution of (\ref{equationz}), with $z_*(0)=z_*(L)=0$. By Taylor expansion at north pole, $$
z_\lambda^{\prime}(0)=-\Phi_\lambda^{\prime \prime}(0)=\frac{\lambda}{N-1}>0 .
$$
If $z_\lambda$ achieve $0$ for the first time at some $T \leq L$, then $z_\lambda>0$ on $(0, T)$. Denote
$$
q=(N-2) \sin ^{N-4} \theta, \, p(\theta):=\sin^{N-2}\theta,
$$
there are two equations:
$$
\begin{aligned}
& -\left(p z_\lambda^{\prime}\right)^{\prime}+q z_\lambda=\lambda p z_\lambda, \\
& -\left(p z_*^{\prime}\right)^{\prime}+q z_*=2 N p z_*,
\end{aligned}
$$
multiply the first one by $z_*$ and the second one by $z_\lambda$ and subtract the two equations, we will have
$$
W^{\prime}:=\left[p\left(z_\lambda z_*^{\prime}-z_* z_\lambda^{\prime}\right)\right]^{\prime}=(\lambda-2 N) p z_\lambda z_*<0
$$ on $(0,T)$. It's also trivial to verify that near $0$, $W\rightarrow0$, so $W<0$ on $(0, T]$. If $T<L$, then $z_\lambda(T)=0, z_\lambda^{\prime}(T)<0$ ( If \(z_\lambda'(T)=0\), then together with \(z_\lambda(T)=0\), uniqueness for the linear initial value problem implies \(z_\lambda\equiv0\), contradicting its positivity on \((0,T)\). ), $z_*(T)>0,$ then $W>0$, which is a contradiction. If $T=L$, then $W=0$, which is also a contradiction. So we have
$$
z_\lambda>0 \quad(0<\theta \leq L), \quad \Phi_\lambda^{\prime}<0 \quad(0<\theta \leq L) .
$$

Then first we choose a $\delta>0$ sufficiently small such that $\Phi_\lambda>\frac{1}{2}$ on $[0, \delta]$, since $\frac{V_a}{a} \rightarrow \Phi_\lambda$, we could choose small $a$ such that
$$
\left|\frac{V_a(\theta)}{a}-\Phi_\lambda(\theta)\right|<\frac{1}{4},
$$
so
$$
V_a(\theta)>\frac{a}{4}>0 ,
$$
then by (\ref{integral}), $V_a^{\prime}<0$. On $[\delta, L]$, since now we could find a $c_\delta$ such that $\Phi_\lambda^{\prime}(\theta) \leq-c_\delta<0, \quad \delta \leq \theta \leq L$, combining with the fact that $\frac{V_a^{\prime}}{a} \rightarrow \Phi_\lambda^{\prime}$ we could finish the proof.
 
\end{proof}

\begin{lemma}
Fix $\lambda>2(N-2)$. For any $a_j \rightarrow \infty$, we could choose a subsequence which tends to $U \in C^2((0, L])$ locally away from the north pole, besides, near north pole
\begin{equation}
U(\theta)=-2 \log \sin \theta+\log \frac{2(N-3)}{\lambda}+\frac{\mu}{4(N-2)} \theta^2+O\left(\theta^4\right).
\end{equation}
\end{lemma}
\begin{proof}
We divide the proof into 5 steps.
\paragraph{\textbf{Step 1: Derivation of the equation for $V_a-V_\lambda^*$.}}

Denote $\mu:=\lambda-2(N-2)>0$ and 
$$
V_\lambda^*(\theta):=-2 \log \sin \theta+\log \frac{2(N-3)}{\lambda},
$$
we have
$$
\Delta_{\mathbb{S}^{N-1}} V_\lambda^*+\lambda\left(e^{V_\lambda^*}-1\right)=-\mu .
$$
Then we do the Emden-Fowler transform:
$$
t=\log \tan \frac{\theta}{2}, \quad-\infty<t \leq 0,
$$
and
$$
x_a(t):=V_a(\theta)-V_\lambda^*(\theta) .
$$
Then the original equation becomes
\begin{equation}\label{transformMAIN}
x_a^{\prime \prime}-(N-3) \tanh t x_a^{\prime}+2(N-3)\left(e^{x_a}-1\right)=\mu \operatorname{sech}^2 t .
\end{equation}

\paragraph{\textbf{Step 2: Rescaling near the north pole and convergence to $u_0$.}}

When $t \rightarrow-\infty$, $\theta \rightarrow 0$, we have
$$
\sin \theta \sim \theta \sim 2 e^t,
$$
then
$$
\log \sin \theta\sim\log \left(2 e^t\right)+o(1)=t+\log 2+o(1),
$$
besides, combining with
$$
V_a(\theta)=a+o(1) \quad(\theta \rightarrow 0),
$$
we will have 
$$
x_a(t)=2 t+a+\gamma_\lambda+o(1), \quad x_a^{\prime}(t) \rightarrow 2, \quad t \rightarrow-\infty,
$$
with 
$$
\gamma_\lambda:=\log \frac{2 \lambda}{N-3} .
$$
Then we set
$$
s=t+\frac{a}{2}, \quad \varepsilon=e^{-a}, \quad r=e^s,
$$
and define
$$
u_{\varepsilon}(r):=x_a\left(t\right)-2 s-\gamma_\lambda .
$$
Then $u_{\varepsilon}(0)=u_{\varepsilon}^{\prime}(0)=0$, and 
\begin{equation}\label{transformODE}
\begin{aligned}
u_{\varepsilon}^{\prime \prime}+\frac{N-2}{r} u_{\varepsilon}^{\prime}+4 \lambda e^{u_{\varepsilon}} & -\frac{2(N-3) \varepsilon}{1+\varepsilon r^2}\left(r u_{\varepsilon}^{\prime}+2\right) \\
& -\frac{4 \mu \varepsilon}{\left(1+\varepsilon r^2\right)^2}=0 .
\end{aligned}
\end{equation}
Define $$
P_{\varepsilon}(r):=r^{N-2}\left(1+\varepsilon r^2\right)^{-(N-3)} .
$$
Then (\ref{transformODE}) becomes
$$
\left(P_{\varepsilon} u_{\varepsilon}^{\prime}\right)^{\prime}=P_{\varepsilon}\left[-4 \lambda e^{u_{\varepsilon}}+C_{\varepsilon}(r)\right],
$$ 
where
$$
C_{\varepsilon}(r)=\frac{4(N-3) \varepsilon}{1+\varepsilon r^2}+\frac{4 \mu \varepsilon}{\left(1+\varepsilon r^2\right)^2} .
$$
So 
$$
u_{\varepsilon}^{\prime}(r)=P_{\varepsilon}(r)^{-1} \int_0^r P_{\varepsilon}(q)\left[-4 \lambda e^{u_{\varepsilon}(q)}+C_{\varepsilon}(q)\right] \mathrm{d} q .
$$
When $r$ is bounded, we could get the limiting initial value problem
\begin{equation}
u_0^{\prime \prime}+\frac{N-2}{r} u_0^{\prime}+4 \lambda e^{u_0}=0, \quad u_0(0)=u_0^{\prime}(0)=0,
\end{equation}
and 
$$
u_0^{\prime}(r)=-4 \lambda r^{-(N-2)} \int_0^r q^{N-2} e^{u_0(q)} \mathrm{d} q<0.
$$
So $u_0 \leq 0$ and $\left|u_0^{\prime}(r)\right| \leq 4 \lambda r /(N-1)$.

Fix $R>0$ and take $\varepsilon>0$ sufficiently small that
$\varepsilon R^2\leq1/2$, so that $u_\varepsilon$ is defined
on $[0,R]$.
Throughout the argument, $C_R$ denotes a positive constant
depending only on $R$, $N$, and $\lambda$, which may change
from line to line but is independent of $\varepsilon$.

Recall the integral identities
\[
u_\varepsilon'(r)
=\int_0^r\frac{P_\varepsilon(q)}{P_\varepsilon(r)}
\left[-4\lambda e^{u_\varepsilon(q)}+C_\varepsilon(q)\right]\,dq,
\,
u_0'(r)
=-4\lambda\int_0^r\left(\frac qr\right)^{N-2}
e^{u_0(q)}\,dq,
\]
where
\[
P_\varepsilon(r)
=r^{N-2}(1+\varepsilon r^2)^{-(N-3)},
\qquad
C_\varepsilon(r)
=\frac{4(N-3)\varepsilon}{1+\varepsilon r^2}
+\frac{4\mu\varepsilon}{(1+\varepsilon r^2)^2}.
\]
For $0\leq q\leq r\leq R$ with $r>0$, we have
\[
\frac{P_\varepsilon(q)}{P_\varepsilon(r)}
=\left(\frac qr\right)^{N-2}
\left(\frac{1+\varepsilon r^2}
{1+\varepsilon q^2}\right)^{N-3}.
\]
Since
\[
1\leq
\frac{1+\varepsilon r^2}{1+\varepsilon q^2}
\leq\frac32,
\qquad
0\leq
\frac{1+\varepsilon r^2}{1+\varepsilon q^2}-1
=\frac{\varepsilon(r^2-q^2)}{1+\varepsilon q^2}
\leq\varepsilon r^2,
\]
the mean value theorem gives
\[
\frac{P_\varepsilon(q)}{P_\varepsilon(r)}
\leq C_R\left(\frac qr\right)^{N-2},
\qquad
\left|
\frac{P_\varepsilon(q)}{P_\varepsilon(r)}
-\left(\frac qr\right)^{N-2}
\right|
\leq C_R\varepsilon r^2
\left(\frac qr\right)^{N-2}.
\]
Moreover,
\[
|C_\varepsilon(q)|\leq C_R\varepsilon
\qquad(0\leq q\leq R).
\]

Set
\[
M(r):=\sup_{0\leq q\leq r}
|u_\varepsilon(q)-u_0(q)|.
\]
Whenever $M(r)\leq1$, the inequality $u_0\leq0$ implies
$u_\varepsilon\leq1$ on $[0,r]$. Hence
\[
|e^{u_\varepsilon(q)}-e^{u_0(q)}|
\leq e\,|u_\varepsilon(q)-u_0(q)|
\leq eM(r)
\qquad(0\leq q\leq r).
\]
Subtracting the two integral identities yields
\[
\begin{aligned}
(u_\varepsilon-u_0)'(r)
={}&-4\lambda\int_0^r
\frac{P_\varepsilon(q)}{P_\varepsilon(r)}
\left(e^{u_\varepsilon(q)}-e^{u_0(q)}\right)\,dq\\
&-4\lambda\int_0^r
\left[
\frac{P_\varepsilon(q)}{P_\varepsilon(r)}
-\left(\frac qr\right)^{N-2}
\right]e^{u_0(q)}\,dq\\
&+\int_0^r
\frac{P_\varepsilon(q)}{P_\varepsilon(r)}
C_\varepsilon(q)\,dq.
\end{aligned}
\]
Using $e^{u_0}\leq1$ and
\[
\int_0^r\left(\frac qr\right)^{N-2}\,dq
=\frac{r}{N-1},
\]
we obtain
\[
|(u_\varepsilon-u_0)'(r)|
\leq C_RrM(r)+C_R\varepsilon r^3+C_R\varepsilon r
\leq C_RrM(r)+C_R\varepsilon r.
\]
The last inequality uses $r\leq R$.
It also holds at $r=0$, since both initial derivatives vanish.

We now show that $M(r)<1$ throughout $[0,R]$.
Since $u_\varepsilon(0)=u_0(0)=0$, if
$|u_\varepsilon-u_0|$ reaches $1$ on $[0,R]$, continuity
gives a first such point $r_*\in(0,R]$ satisfying
\[
|u_\varepsilon(r_*)-u_0(r_*)|=1,
\qquad M(r)\leq1\quad(0\leq r\leq r_*).
\]
The derivative estimate therefore applies on $[0,r_*]$.
Integrating from $0$ to any $s\leq r_*$ gives
\[
|u_\varepsilon(s)-u_0(s)|
\leq\frac{C_R\varepsilon s^2}{2}
+C_R\int_0^s qM(q)\,dq.
\]
Taking the supremum over $0\leq s\leq r$ yields
\[
M(r)\leq\frac{C_RR^2}{2}\varepsilon
+C_R\int_0^r qM(q)\,dq,
\qquad 0\leq r\leq r_*.
\]
Fixing a constant $C_R$ for which this inequality holds,
Grönwall's inequality gives
\[
M(r)\leq
\frac{C_RR^2}{2}\varepsilon
\exp\left(\frac{C_Rr^2}{2}\right)
\leq K_R\varepsilon,
\qquad
K_R:=\frac{C_RR^2}{2}
\exp\left(\frac{C_RR^2}{2}\right).
\]
In particular,
\[
1=|u_\varepsilon(r_*)-u_0(r_*)|
\leq M(r_*)\leq K_R\varepsilon<1
\]
for sufficiently small $\varepsilon$, a contradiction.
Here $K_R$ is independent of $\varepsilon$ and $r_*$.

Consequently, $M(r)<1$ on $[0,R]$.
The same derivative estimate and Grönwall argument
now apply on the whole interval, giving
\[
M(R)\leq K_R\varepsilon.
\]
Substituting this bound into the derivative estimate yields
\[
\sup_{0\leq r\leq R}|u_\varepsilon'(r)-u_0'(r)|
\leq C_R\varepsilon.
\]
Therefore,
\[
\|u_\varepsilon-u_0\|_{C^1([0,R])}
\leq C_R\varepsilon\longrightarrow0
\qquad(\varepsilon\downarrow0).
\]

\paragraph{\textbf{Step 3: Long-time behavior of the autonomous orbit.}}

Define
$$
\bar{x}(s):=u_0\left(e^s\right)+2 s+\gamma_\lambda ,
$$
then if we set
\begin{equation}
\bar{x}^{\prime \prime}+(N-3) \bar{x}^{\prime}+2(N-3)\left(e^{\bar{x}}-1\right)=0 ,
\end{equation}
and
\begin{equation}
\mathcal{H}(x, y):=\frac{1}{2} y^2+2(N-3)\left(e^x-x-1\right) ,
\end{equation}
then
\begin{equation}
\frac{\mathrm{d}}{\mathrm{~d} s} \mathcal{H}\left(\bar{x}, \bar{x}^{\prime}\right)=-(N-3)\left(\bar{x}^{\prime}\right)^2 \leq 0 .
\end{equation}
Since $\mathcal{H} \geq 0$, there exists $E_{\infty} \geq 0$ such that
$$
\mathcal{H}\left(\bar{x}(s), \bar{x}^{\prime}(s)\right) \rightarrow E_{\infty} .
$$
Besides, since $e^x-x-1 \rightarrow \infty$ as $x \rightarrow \pm \infty$,
we have the boundness of $\bar{x}, \bar{x}^{\prime}$, by the equation itself, we have the boundness of $\bar{x}^{\prime \prime}$. Combining with the fact that
$$
(N-3) \int_{s_0}^{\infty}\left(\bar{x}^{\prime}\right)^2 d s=\mathcal{H}\left(\bar{x}\left(s_0\right), \bar{x}^{\prime}\left(s_0\right)\right)-E_{\infty}<\infty,
$$
we will have 
$$
\bar{x}^{\prime}(s) \rightarrow 0 .
$$
Since $\bar{x}'(s)\to0$ and
$\mathcal{H}(\bar{x}(s),\bar{x}'(s))\to E_\infty$, we have
\[
e^{\bar{x}(s)}-\bar{x}(s)-1
\longrightarrow \frac{E_\infty}{2(N-3)}.
\]
If $E_\infty=0$, this immediately implies $\bar{x}(s)\to0$.
If $E_\infty>0$, then $\bar{x}(s)$ stays away from zero
for all sufficiently large $s$ and, by continuity, has a
fixed sign. Since $x\mapsto e^x-x-1$ is strictly monotone
on each of the intervals $(-\infty,0)$ and $(0,\infty)$,
it follows that $\bar{x}(s)\to\ell$ for some $\ell\neq0$.

Then it's trivial to prove that if $\bar{x}(s) \rightarrow \ell \neq 0$, then $$
\bar{x}^{\prime \prime}(s)=-(N-3) \bar{x}^{\prime}(s)-2(N-3)\left(e^{\bar{x}(s)}-1\right) \longrightarrow-2(N-3)\left(e^{\ell}-1\right) \neq 0,
$$
which is impossible, so
\begin{equation}
\left(\bar{x}(s), \bar{x}^{\prime}(s)\right) \longrightarrow(0,0), \quad s \rightarrow+\infty .
\end{equation}

\paragraph{\textbf{Step 4: Construction of the singular limit.}}

Recall that
$$
x_a^{\prime \prime}-(N-3) \tanh t x_a^{\prime}+2(N-3)\left(e^{x_a}-1\right)=\mu \operatorname{sech}^2 t ,
$$
we have
\begin{equation}
\frac{\mathrm{d}}{\mathrm{~d} t} \mathcal{H}\left(x_a, x_a^{\prime}\right)=(N-3) \tanh t\left(x_a^{\prime}\right)^2+\mu \operatorname{sech}^2 t x_a^{\prime} .
\end{equation}
Since $t \leq 0$, $\tanh t \leq 0$, besides, $\left|x_a^{\prime}\right| \leq \sqrt{2 \mathcal{H}}$, so
\begin{equation}\label{inte}
\frac{\mathrm{d}}{\mathrm{~d} t} \sqrt{\mathcal{H}\left(x_a, x_a^{\prime}\right)} \leq \frac{\mu}{\sqrt{2}} \operatorname{sech}^2 t .
\end{equation}
Now we fix $T<0$ and $S>0$, then set 
$$
t_a:=S-\frac{a}{2} .
$$
Recall the definition
$$
s=t+\frac{a}{2}, \quad \varepsilon=e^{-a}, \quad r=e^s, u_{\varepsilon}(r)=x_a\left(s-\frac{a}{2}\right)-2 s-\gamma_\lambda .
$$
Now if we fix $S$, then
$$
s=t_a+\frac{a}{2}=S, u_{\varepsilon}\left(e^S\right)=x_a\left(S-\frac{a}{2}\right)-2 S-\gamma_\lambda .
$$
Then
$$
u_{\varepsilon}\left(e^S\right)=x_a\left(S-\frac{a}{2}\right)-2 S-\gamma_\lambda, x_a\left(t_a\right)=u_{\varepsilon}\left(e^S\right)+2 S+\gamma_\lambda .
$$

Recall that
$$
\bar{x}(S)=u_0\left(e^S\right)+2 S+\gamma_\lambda ,
$$
then when $a \rightarrow \infty$ 
\begin{equation}\label{ainfty}
\left(x_a\left(t_a\right), x_a^{\prime}\left(t_a\right)\right) \rightarrow\left(\bar{x}(S), \bar{x}^{\prime}(S)\right) .
\end{equation}
integrate (\ref{inte}) from $t_a$ to $T$, if $t_a=S-a / 2<T$, we will have
\begin{equation}
\begin{aligned}
\sqrt{H\left(x_a(T), x_a^{\prime}(T)\right)} \leq & \sqrt{H\left(x_a(t_a), x_a^{\prime}(t_a)\right)} \\
& +\frac{\mu}{\sqrt{2}}[\tanh T-\tanh (t_a)] \\
\leq & \sqrt{H\left(x_a(t_a), x_a^{\prime}(t_a)\right)}+\frac{\mu}{\sqrt{2}}(1+\tanh T) .
\end{aligned}
\end{equation}
If $H$ vanishes, we follow the same argument above with
$\sqrt{H+\eta}$ in place of $\sqrt H$, where $\eta>0$,
and let $\eta\downarrow0$ after integration to obtain
the same estimates.

For any $\delta$, we could find $S$ such that
$$
\sqrt{H\left(\bar{x}(S), \bar{x}^{\prime}(S)\right)}<\delta,
$$
then by (\ref{ainfty}), there exists $A=A(S, \delta, T)$ such that for all $a>A$ we have
$$
S-a / 2<T, \quad \sqrt{H\left(x_a(t_a), x_a^{\prime}(t_a)\right)}<2 \delta .
$$
So
$$
\sqrt{H\left(x_a(T), x_a^{\prime}(T)\right)} \leq 2 \delta+\frac{\mu}{\sqrt{2}}(1+\tanh T), \quad a>A .
$$
By the arbitrariness of $\delta$, we will have 
\begin{equation}\label{random}
\limsup _{a \rightarrow \infty} \sqrt{H\left(x_a(T), x_a^{\prime}(T)\right)} \leq \frac{\mu}{\sqrt{2}}(1+\tanh T) .
\end{equation}
For each fixed $T<0$, the functions $x_a$, $x_a'$ and
$x_a''$ are uniformly bounded on $[T,0]$ for sufficiently
large $a$. In particular, $x_a$ and $x_a'$ are equicontinuous.
Given any sequence $a_j\to+\infty$, the Arzel\`a--Ascoli
theorem and a diagonal extraction on $[-k,0]$, $k\in\mathbb N$,
yield a subsequence, still denoted by $a_j$, such that
\[
x_{a_j}\longrightarrow x_\infty
\quad\text{in }C^1_{\mathrm{loc}}((-\infty,0]).
\]
Moreover, the equation gives
\[
x_{a_j}''
=(N-3)\tanh t\,x_{a_j}'
-2(N-3)(e^{x_{a_j}}-1)+\mu\operatorname{sech}^2t.
\]
The right-hand side converges uniformly on every compact
subinterval of $(-\infty,0]$. Hence $x_\infty$ is $C^2$ and
\[
x_{a_j}\longrightarrow x_\infty
\quad\text{in }C^2_{\mathrm{loc}}((-\infty,0]).
\]
The limiting function satisfies
\begin{equation}\label{limitfunction}
x_{\infty}^{\prime \prime}-(N-3) \tanh t x_{\infty}^{\prime}+2(N-3)\left(e^{x_{\infty}}-1\right)=\mu \operatorname{sech}^2 t .
\end{equation}

\paragraph{\textbf{Step 5: Asymptotic expansion near the north pole.}}

By (\ref{random}) and the definition of $\mathcal{H}$, we have
\begin{equation}\label{t_to_infty}
x_{\infty}(t)=O\left(e^{2 t}\right), \quad x_{\infty}^{\prime}(t)=O\left(e^{2 t}\right), \quad t \rightarrow-\infty .
\end{equation}
(\ref{limitfunction}) could be written as
$$
L x_{\infty}=4 \mu e^{2 t}+O\left(e^{4 t}\right), \quad L:=D^2+(N-3) D+2(N-3) .
$$
We define
$$
z(t):=x_{\infty}(t)-\frac{\mu}{N-2} e^{2 t},
$$
then
\begin{equation}\label{approximate_equation}
z^{\prime \prime}+(N-3) z^{\prime}+2(N-3) z=f(t), \quad|f(t)| \leq C e^{4 t} .
\end{equation}
Write (\ref{approximate_equation}) as a system, we get
$$
Z(t):=\binom{z(t)}{z^{\prime}(t)}, \quad Z^{\prime}=A Z+\binom{0}{f(t)}, \quad A=\left(\begin{array}{cc}
0 & 1 \\
-2(N-3) & -(N-3)
\end{array}\right) .
$$
Since all eigenvalues of \(A\) have strictly negative real parts, its Jordan normal form implies that there exists $M, \kappa>0$, such that
$$
\left\|e^{A u}\right\| \leq M e^{-\kappa u}, \quad u \geq 0 .
$$
Variation of constant gives
$$
Z(t)=e^{A(t-s)} Z(s)+\int_s^t e^{A(t-r)}\binom{0}{f(r)} d r,
$$
if we fix $t$ and let $s \rightarrow -\infty$, we get
$$
\left\|e^{A(t-s)} Z(s)\right\| \leq M C e^{-\kappa(t-s)} e^{2 s}=M C e^{-\kappa t} e^{(\kappa+2) s} \longrightarrow 0,
$$
then
$$
Z(t)=\int_{-\infty}^t e^{A(t-r)}\binom{0}{f(r)} d r .
$$
Combining with the fact that $|f(r)| \leq C e^{4 r}$ we will have
$$
z(t)=O\left(e^{4 t}\right), \quad z^{\prime}(t)=O\left(e^{4 t}\right),
$$
then
$$
x_{\infty}(t)=\frac{\mu}{N-2} e^{2 t}+O\left(e^{4 t}\right) .
$$

The equation for $z$, together with
$z(t),z'(t),f(t)=O(e^{4t})$, also gives
$z''(t)=O(e^{4t})$ as $t\to-\infty$.
Returning to the original coordinate, define
\[
t(\theta):=\log\tan\frac{\theta}{2},
\qquad
w(\theta):=x_\infty(t(\theta)),
\qquad
U(\theta):=V_\lambda^*(\theta)+w(\theta).
\]
Since the coordinate change is smooth on every compact
subinterval of $(0,L]$, the convergence of $x_{a_j}$ implies
\[
V_{a_j}\longrightarrow U
\quad\text{in }C^2_{\mathrm{loc}}((0,L]).
\]
In particular, $U$ satisfies the original equation on $(0,L]$.

By the definition of $z$,
\[
w(\theta)
=\frac{\mu}{N-2}\tan^2\frac{\theta}{2}
+z(t(\theta)).
\]
Moreover,
\[
t'(\theta)=\frac{1}{\sin\theta},
\qquad
t''(\theta)=-\frac{\cos\theta}{\sin^2\theta},
\qquad
e^{4t(\theta)}=O(\theta^4).
\]
Thus the estimates for $z$, $z'$, and $z''$, together with
the chain rule, yield
\[
z(t(\theta))=O(\theta^4),\qquad
\frac{d}{d\theta}z(t(\theta))=O(\theta^3),\qquad
\frac{d^2}{d\theta^2}z(t(\theta))=O(\theta^2).
\]
Consequently, as $\theta\downarrow0$,
\[
\begin{aligned}
w(\theta)
&=\frac{\mu}{4(N-2)}\theta^2+O(\theta^4),\\
w'(\theta)
&=\frac{\mu}{2(N-2)}\theta+O(\theta^3),\\
w''(\theta)
&=\frac{\mu}{2(N-2)}+O(\theta^2).
\end{aligned}
\]
Substituting the expansion of $w$ into
$U=V_\lambda^*+w$ proves the asserted asymptotic expansion.

\end{proof}

\begin{lemma}
If $N \geq 4$ and $\lambda>2(N-2)$ then any singular limit $U$ has strict minimum point at some $T \in(0, \pi / 2)$.
\end{lemma}
\begin{proof}
By Lemma 4.2 and the derivative estimates above,
\[
U'(\theta)=-2\cot\theta+w'(\theta)
=-\frac{2}{\theta}+O(\theta)<0
\]
for all sufficiently small $\theta>0$.

If $U'$ has a first zero $T\in(0,L)$, then $U''(T)\ge0$.
If $U''(T)=0$, the equation implies $U(T)=0$.
Together with $U'(T)=0$, uniqueness for the regular initial
value problem at $T$ would then give $U\equiv0$ on $(0,L)$,
contradicting its singular behavior at the north pole.
Thus $U''(T)>0$, and $T$ is a strict local minimum of $U$.
It therefore suffices to rule out the possibility that
$U'<0$ throughout $(0,L)$.

We assume by contradiction that
$$
U^{\prime}(\theta)<0, \quad 0<\theta<L .
$$
Besides, we set
\begin{equation}
\phi(\theta):=-V_\lambda^{* \prime}(\theta)=2 \cot \theta>0, \quad v:=-U^{\prime}>0, \quad \rho:=\frac{v}{\phi} ,
\end{equation}
since $U=V_\lambda^*+w$, we have
$$
w^{\prime}=\phi(1-\rho) .
$$
Combining with the fact that near north pole:
\begin{equation}
U(\theta)=V_\lambda^*(\theta)+w(\theta), \quad w(\theta)=\frac{\mu}{4(N-2)} \theta^2+O\left(\theta^4\right) ,
\end{equation}
we will have that near north pole:
\begin{equation}
\rho(\theta)=1-\frac{\mu}{4(N-2)} \theta^2+O\left(\theta^4\right),
\end{equation}
so near north pole, we have
\begin{equation}
w>0, \quad 0<\rho<1 .
\end{equation}
Set $v=-U^{\prime}$, then 
\begin{equation}
v^{\prime \prime}+(N-2) \cot \theta v^{\prime}+\left[\frac{2(N-3) e^w}{\sin ^2 \theta}-(N-2) \csc ^2 \theta\right] v=0 .
\end{equation}
Set $\phi=-V_\lambda^{* \prime}$, then
\begin{equation}
\phi^{\prime \prime}+(N-2) \cot \theta \phi^{\prime}+\left[\frac{2(N-3)}{\sin ^2 \theta}-(N-2) \csc ^2 \theta\right] \phi=0 .
\end{equation}
Then we have
\begin{equation}
\left(\sin ^{N-2} \theta \phi^2 \rho^{\prime}\right)^{\prime}=-2(N-3) \sin ^{N-4} \theta \phi^2\left(e^w-1\right) \rho .
\end{equation}
Then by $\rho(\theta)=1-\frac{\mu}{4(N-2)} \theta^2+O\left(\theta^4\right)$, then near north pole:
\begin{equation}
\sin ^{N-2} \theta \phi^2 \rho^{\prime}=O\left(\theta^{N-3}\right) \rightarrow 0 \quad(\theta \downarrow 0).
\end{equation}
On any region, if $w>0$ and $\rho>0$, then
$$
\left(\sin ^{N-2} \theta \phi^2 \rho^{\prime}\right)^{\prime}<0.
$$

If $w$ had a first zero $\theta_*\in(0,L)$, this argument
would give $w'>0$ on $(0,\theta_*)$, so that
$w(\theta_*)\geq w(\theta_1)>0$ for any
$\theta_1\in(0,\theta_*)$, a contradiction.
Consequently,
\[
w>0,\qquad 0<\rho<1,\qquad \rho'<0
\quad\text{throughout }(0,L).
\]

For fixed $\theta_0 \in(0, L)$, there exists $c_0>0$ such that 
\begin{equation}
\sin ^{N-2} \theta \phi^2 \rho^{\prime}(\theta) \leq-c_0, \quad \theta \in\left[\theta_0, L\right) .
\end{equation}
When $\theta \uparrow L$, $$
\phi(\theta)=2 \cot \theta \sim 2(L-\theta),
$$
$$
\sin ^{N-2} \theta \phi^2 \sim 4(L-\theta)^2 .
$$
So around $L$ we have
\begin{equation}
\rho^{\prime}(\theta) \leq-\frac{c}{(L-\theta)^2}.
\end{equation}
So after integrating, $\rho$ will become negative before $L$, which contradicts to $\rho=v / \phi>0$. So there exists the first $T<L$ such that 
$$
U^{\prime}(T)=0 .
$$
So $U^{\prime \prime}(T) \geq 0$, uniqueness for the regular initial value problem at \(T\) then forces 
\[
U''(T)>0.
\]
\end{proof}

\begin{lemma}
If $N \geq 4$ and $2(N-2)<\lambda<2 N$ then there exists $a_* \in(0, \infty)$ such that $V_{a_*}$ is smooth on $[0, L]$ and
$$
V_{a_*}^{\prime}(0)=V_{a_*}^{\prime}(L)=0 .
$$
\end{lemma}
\begin{proof}
From last lemma, we could choose
$$
0<\theta_{-}<T<\theta_{+}<L
$$
such that 
$$
U^{\prime}\left(\theta_{-}\right)<0, \quad U^{\prime}\left(\theta_{+}\right)>0 .
$$
From Lemma 4.2 we could have 
$$
V_{a_j}^{\prime}\left(\theta_{-}\right)<0, \quad V_{a_j}^{\prime}\left(\theta_{+}\right)>0 ,
$$
for sufficiently large $a_j$.

Define
$$
\mathcal{S}:=\left\{a>0: V_a^{\prime}(\theta)<0 \text { for every } \theta \in(0, L]\right\}, \, a_*:=\sup \{A>0:(0, A) \subset \mathcal{S}\} .
$$

Lemma 4.1 ensures that $a_*>0$, while Lemmas 4.2 and 4.3 provide a positive initial value outside $\mathcal{S}$, so $a_*<\infty$. Taking $a_j \uparrow a_*$ with $a_j<a_*$, continuous dependence on the initial value yields
$$
V_{a_*}^{\prime}(\theta) \leq 0 \quad(0<\theta \leq L) .
$$

If there exists $\theta_0 \in(0, L)$ such that $V_{a_*}^{\prime}\left(\theta_0\right)=0$, then 
$$
V_{a_*}^{\prime \prime}\left(\theta_0\right)=0 .
$$
Then from the original equation we have
$$
V_{a_*}\left(\theta_0\right)=V_{a_*}^{\prime}\left(\theta_0\right)=0.
$$
Then $V_{a_*} \equiv 0$, a contradiction. So
\begin{equation}
V_{a_*}^{\prime}(\theta)<0, \quad 0<\theta<L .
\end{equation}
Since $V_{a_{*}}(0)=a_*>0$, there exists $\delta$ such that
\begin{equation}
V_{a_{*}}(\theta)>\frac{3 a_*}{4}, \quad 0 \leq \theta \leq \delta .
\end{equation}
\textbf{If we assume $V^{\prime}_{a_{*}}(L)<0$.} Then by choosing $$
-\eta:=\max _{\theta \in[\delta, L]} V_{a *}^{\prime}(\theta) ,
$$
we will have
$$
V_{a_{*}}^{\prime}(\theta) \leq -\eta, \quad \delta \leq \theta \leq L ,
$$

Then since when $a \rightarrow a_*$, we have
$$
V_a \rightarrow V_{a_*} \quad \text { in } C^1([0, L]).
$$
Therefore, there exists $\varepsilon$ such that,
for every $a$ satisfying $|a-a_*|<\varepsilon$, we have
$$
\sup _{0 \leq \theta \leq \delta}\left|V_a(\theta)-V_{a_*}(\theta)\right|<\frac{a_*}{4},
$$
and
$$
\sup _{\delta \leq \theta \leq L}\left|V_a^{\prime}(\theta)-V_{a_*}^{\prime}(\theta)\right|<\frac{\eta}{2} ,
$$
recalling the fact that $$
V_{a_*}(\theta)>\frac{3 a_*}{4}, \quad 0 \leq \theta \leq \delta ,
$$
we have
$$
V_a(\theta) \geq V_{a_*}(\theta)-\left|V_a(\theta)-V_{a_*}(\theta)\right|>\frac{3 a_*}{4}-\frac{a_*}{4}=\frac{a_*}{2}>0 .
$$
So we have
$$
V_a(\theta)>0, \quad 0 \leq \theta \leq \delta .
$$
On the other hand,
$$
\begin{aligned}
V_a^{\prime}(\theta) & \leq V_{a_*}^{\prime}(\theta)+\left|V_a^{\prime}(\theta)-V_{a_*}^{\prime}(\theta)\right| \\
& \leq-\eta+\frac{\eta}{2}=-\frac{\eta}{2}<0 ,
\end{aligned}
$$
so we have
$$
V_a^{\prime}(\theta)<0, \quad \delta \leq \theta \leq L .
$$
By the structure 
$$
V_a^{\prime}(\theta)=-\lambda \sin ^{-(N-2)} \theta \int_0^\theta \sin ^{N-2} s\left(e^{V_a(s)}-1\right) d s,
$$
we will also have 
$$
V_a^{\prime}(\theta)<0, \quad 0<\theta \leq \delta .
$$
We will have
$$
V_a^{\prime}(\theta)<0, \quad 0<\theta \leq L .
$$
Thus every $a$ with $|a-a_*|<\varepsilon$ belongs to $S$.
Together with $(0,a_*)\subset S$, this gives
\[
(0,a_*+\varepsilon)\subset S,
\]
contradicting the definition of $a_*$.
Therefore $V_{a_*}'(L)=0$.
\end{proof}

\begin{remark}
A combination of the mass identity, spectral and flux
estimates, and rescaling near the poles with
Emden--Fowler energy estimates may lead to a
nonexistence result for nonconstant smooth solutions
when $\lambda>0$ is sufficiently small, within the
axially symmetric and equatorially symmetric class.
However, such a result would serve here only as
qualitative context for the numerical diagrams.
To establish existence, it suffices to construct
a nonconstant solution within a symmetry class.
A general nonexistence result must also rule out
nonconstant solutions outside that class.
Excluding symmetric solutions alone would therefore
leave the general nonexistence problem unresolved.

\end{remark}

\section{Pin initial value $V(0)=a<0$}
\begin{theorem}\label{thm:existence-above-2N}
Let $N\ge4$ and $\lambda>2N$. Then
\[
\Delta_{S^{N-1}}V+\lambda(e^V-1)=0
\quad\text{on }S^{N-1}
\]
admits a nonconstant smooth solution which is axially symmetric
and symmetric about the equator. Moreover, the solution can be
chosen so that
\[
V(0)<0,
\qquad
V'(\theta)>0
\quad\text{for }0<\theta<\frac{\pi}{2}.
\]
\end{theorem}

\begin{proof}
Fix $\lambda>2N$ and set $L=\pi/2$. For each $a<0$, let $V_a$
be the regular solution of
\[
V_a''+(N-2)\cot\theta\,V_a'
+\lambda(e^{V_a}-1)=0,
\qquad
V_a(0)=a,\quad V_a'(0)=0.
\]
The energy
\[
E_a(\theta)
:=\frac12(V_a'(\theta))^2
+\lambda\bigl(e^{V_a(\theta)}-V_a(\theta)-1\bigr)
\]
satisfies
\[
E_a'(\theta)=-(N-2)\cot\theta\,(V_a'(\theta))^2\le0
\quad(0<\theta<L).
\]
Consequently, every finite initial value produces a solution
on $[0,L]$. The integral formulation also gives continuous
dependence on $a$ in $C^1([0,L])$.

We first show that sufficiently negative initial values give
solutions which remain strictly increasing up to the equator.
Set
\[
p(\theta):=\sin^{N-2}\theta.
\]
The integral identity
\begin{equation}\label{eq:negative-shooting-integral}
V_a'(\theta)
=\frac{\lambda}{p(\theta)}
\int_0^\theta p(s)\bigl(1-e^{V_a(s)}\bigr)\,ds
\end{equation}
implies
\[
V_a'(\theta)
\le\frac{\lambda}{p(\theta)}\int_0^\theta p(s)\,ds.
\]
Integrating once more yields
\[
V_a(\theta)\le a+\lambda K_N
\quad(0\le\theta\le L),
\]
where
\[
K_N:=
\int_0^L
\frac{\int_0^t\sin^{N-2}s\,ds}{\sin^{N-2}t}\,dt
<\infty.
\]
Indeed, the integrand is asymptotic to $t/(N-1)$ as
$t\downarrow0$. Thus, if $a<-\lambda K_N$, then
$V_a<0$ on $[0,L]$, and
\eqref{eq:negative-shooting-integral} gives
\begin{equation}\label{eq:large-negative-increasing}
V_a'(\theta)>0
\quad(0<\theta\le L).
\end{equation}

We next show that negative initial values sufficiently close
to zero produce a turning point before the equator.
By linearization at the zero solution,
\begin{equation}\label{eq:negative-shooting-linearization}
\frac{V_a}{a}\longrightarrow\Phi_\lambda
\quad\text{in }C^1([0,L])
\quad\text{as }a\uparrow0,
\end{equation}
where
\[
\Phi_\lambda''+(N-2)\cot\theta\,\Phi_\lambda'
+\lambda\Phi_\lambda=0,
\qquad
\Phi_\lambda(0)=1,\quad\Phi_\lambda'(0)=0.
\]
The convergence in \eqref{eq:negative-shooting-linearization}
follows from the integral formulation, the estimate
$\|V_a\|_{C^1([0,L])}=O(|a|)$, and Gr\"onwall's inequality.

Define $z_\lambda:=-\Phi_\lambda'$. Then
\[
-(pz_\lambda')'
+(N-2)\sin^{N-4}\theta\,z_\lambda
=\lambda pz_\lambda,
\]
and
\[
z_\lambda(0)=0,\qquad
z_\lambda'(0)=\frac{\lambda}{N-1}>0.
\]
The function
\[
z_*(\theta):=\sin\theta\cos\theta
\]
satisfies the same equation with $\lambda$ replaced by $2N$.
Hence
\[
W(\theta):=
p(\theta)\bigl(z_\lambda z_*'-z_*z_\lambda'\bigr)
\]
satisfies
\[
W'=(\lambda-2N)pz_\lambda z_*,
\qquad
\lim_{\theta\downarrow0}W(\theta)=0.
\]

Suppose that $z_\lambda$ has no zero in $(0,L)$.
Then $z_\lambda>0$ on $(0,L)$, so $\lambda>2N$ implies
\[
W(L)
=(\lambda-2N)\int_0^L pz_\lambda z_*\,d\theta>0.
\]
On the other hand,
\[
p(L)=1,\qquad z_*(L)=0,\qquad z_*'(L)=-1,
\]
and therefore
\[
W(L)=-z_\lambda(L)\le0,
\]
a contradiction. Thus $z_\lambda$ has a first zero
$T\in(0,L)$. Its positivity on $(0,T)$ and uniqueness
for the linear initial value problem imply
\[
z_\lambda'(T)<0.
\]
We may therefore choose $\theta_+\in(T,L)$ such that
\[
\Phi_\lambda'(\theta_+)=-z_\lambda(\theta_+)>0.
\]
By \eqref{eq:negative-shooting-linearization}, for all
sufficiently small negative $a$,
\[
\frac{V_a'(\theta_+)}{a}>0.
\]
Since $a<0$, this gives
\begin{equation}\label{eq:small-negative-turning}
V_a'(\theta_+)<0.
\end{equation}
In contrast, $V_a'>0$ near the north pole, since
$V_a(0)=a<0$ and \eqref{eq:negative-shooting-integral} applies.

Define
\[
\mathcal S_-:=
\left\{
a<0:
V_a'(\theta)>0
\text{ for every }\theta\in(0,L]
\right\},
\]
and
\[
a_*:=
\sup\left\{
A<0:(-\infty,A)\subset\mathcal S_-
\right\}.
\]
By \eqref{eq:large-negative-increasing}, all sufficiently
negative initial values belong to $\mathcal S_-$.
By \eqref{eq:small-negative-turning}, all negative initial
values sufficiently close to zero lie outside $\mathcal S_-$.
Consequently,
\[
-\infty<a_*<0.
\]

Taking $a_j\uparrow a_*$ with $a_j<a_*$ and using
continuous dependence on the initial value, we obtain
\begin{equation}\label{eq:negative-threshold-nondecreasing}
V_{a_*}'(\theta)\ge0
\quad(0<\theta\le L).
\end{equation}
If $V_{a_*}'(\theta_0)=0$ for some $\theta_0\in(0,L)$,
then \eqref{eq:negative-threshold-nondecreasing} shows that
$\theta_0$ is a local minimum of $V_{a_*}'$, so
$V_{a_*}''(\theta_0)=0$. The equation then gives
$V_{a_*}(\theta_0)=0$. Uniqueness for the regular initial
value problem at $\theta_0$ yields $V_{a_*}\equiv0$,
contradicting $V_{a_*}(0)=a_*<0$. Thus
\begin{equation}\label{eq:negative-threshold-strictly-increasing}
V_{a_*}'(\theta)>0
\quad(0<\theta<L).
\end{equation}

Suppose now that $V_{a_*}'(L)>0$.
Choose $\delta>0$ sufficiently small that
$V_{a_*}<0$ on $[0,\delta]$.
By \eqref{eq:negative-threshold-strictly-increasing}
and compactness,
\[
\min_{\theta\in[\delta,L]}V_{a_*}'(\theta)>0.
\]
Continuous dependence implies that, for all $a$ sufficiently
close to $a_*$,
\[
V_a<0\quad\text{on }[0,\delta],
\qquad
V_a'>0\quad\text{on }[\delta,L].
\]
Equation \eqref{eq:negative-shooting-integral} then gives
$V_a'>0$ also on $(0,\delta]$. Hence a whole neighborhood
of $a_*$, chosen to lie in $(-\infty,0)$, is contained
in $\mathcal S_-$. Together with
$(-\infty,a_*)\subset\mathcal S_-$, this contradicts
the definition of $a_*$.
Combining this with \eqref{eq:negative-threshold-nondecreasing},
we conclude that
\[
V_{a_*}'(L)=0.
\]

Finally, extend $V_{a_*}$ to $[0,\pi]$ by
\[
V(\theta)=
\begin{cases}
V_{a_*}(\theta),&0\le\theta\le L,\\
V_{a_*}(\pi-\theta),&L\le\theta\le\pi.
\end{cases}
\]
The equation is invariant under $\theta\mapsto\pi-\theta$.
Since $V_{a_*}'(L)=0$, uniqueness at the equator ensures
that this extension is smooth there.
Regularity at the north pole, as established at the beginning
of Section 4, and at the south pole by reflection gives a
smooth axially symmetric solution on $S^{N-1}$. It is nonconstant by
\eqref{eq:negative-threshold-strictly-increasing},
and satisfies $V(0)=a_*<0$.
\end{proof}

\medskip
\noindent
\textit{\textbf{Proof of Theorem~\ref{thm:nonconstant-symmetric-solutions}.}}
For $2(N-2)<\lambda<2N$, the construction in
Lemma~4.4 yields a regular solution with positive
pole value, strictly decreasing on $(0,\pi/2)$,
and satisfying $V'(0)=V'(\pi/2)=0$.
Reflection across the equator, together with the pole
regularity established at the beginning of Section 4,
gives the required nonconstant smooth solution on $S^{N-1}$.
\hfill\qedsymbol

\section{Conclusion}

Motivated by prescribed Chern scalar curvature problems,
we have studied the equation
\[
\Delta_{\mathbb S^{N-1}}V+\lambda(e^V-1)=0
\]
with $\lambda>0$ on the high dimensional sphere.
For $N\ge4$ and
$\lambda\in(2(N-2),\infty)\setminus\{2N\}$,
we establish the existence of a nonconstant smooth
solution which is axially symmetric and symmetric
about the equator. The solution can be chosen with
a positive pole value and strict decrease on
$(0,\pi/2)$ when $2(N-2)<\lambda<2N$, and with
a negative pole value and strict increase on this
interval when $\lambda>2N$.
The proof uses the pole value as a shooting parameter
and combines linearization, singular-limit analysis,
and continuous dependence of the shooting solutions.
The numerical computations complement these results
by illustrating the different behaviors associated
with positive and negative pole values.

The construction relies on the symmetry of the sphere. Axial symmetry reduces the equation to an ODE,
while equatorial reflection allows a solution on a
hemisphere to be extended to the whole sphere.
Consequently, the present analysis selects a particular
class of solutions and does not describe the structure
of non-axisymmetric solutions. Moreover, the shooting
and reflection arguments cannot be transferred directly
to manifolds lacking the corresponding geometric
symmetries. Extending the analysis to general Hermitian
backgrounds would also require additional control of
the first-order drift term in the Chern Laplacian.

Further work should therefore investigate
non-axisymmetric solutions. Our purpose here is to explore the positive-parameter
problem on a simple geometric model, where the
dependence on the parameter and the behavior of
shooting solutions can be examined in detail.
The sphere provides a setting in which analytical
results can be compared with numerical observations,
helping to formulate questions about solution
profiles and bifurcation behavior.
We view the present results as a step toward
understanding the positive-parameter regime in
broader geometric settings. We hope that the
existence results established here, together with
the phenomena suggested by the numerical computations,
will help guide the development of methods
applicable beyond the imposed symmetries.

\section*{Acknowledgements}
The authors thank the anonymous referees for their careful
reading and valuable suggestions.

\section*{Declarations}

\noindent\textbf{Funding.} The research is partially supported
by Research Development Fund at XJTLU with Nos. RDF-TP-0014 and
RDS10120260083.

\medskip
\noindent\textbf{Competing interests.} The author declares no
competing interests.

\medskip
\noindent\textbf{Data availability.} No data was used for the
research described in this article.

\medskip
\noindent\textbf{Declaration of generative AI and AI-assisted
technologies in the manuscript preparation process.} The initial
idea and framework of the proof were provided by the authors.
The AI models used were GPT's 5.6 Sol and 6.0 Astra, which
contributed to some technical details of the proof. The authors
have reviewed the entire process and takes full responsibility
for the content.

\end{document}